\documentclass[a4paper,12pt]{article}
\usepackage{a4wide}
\usepackage{amsmath}
\usepackage{amssymb}
\usepackage{amsthm}
\usepackage{latexsym}
\usepackage{graphicx}
\usepackage[english]{babel}
\usepackage{makeidx}

\newtheorem{obs} [subsection]{Remark}

\newtheorem{prop}[subsection]{Proposition}

\newtheorem{teor}[subsection]{Theorem}
\newtheorem{lema}[subsection]{Lemma}
\newtheorem{cor} [subsection]{Corollary}

\begin{document}

\selectlanguage{english} \frenchspacing

\large
\begin{center}
\textbf{A note on the generic initial ideal for complete
intersections.}

\normalsize by \large

Mircea Cimpoea\c s
\end{center}

\normalsize

\footnotetext[1]{This paper was supported by the CEEX Program of the
Romanian Ministry of Education and Research, Contract
CEX05-D11-11/2005 and by the Higher Education Commission of
Pakistan.}

\begin{abstract}
We prove that the $d$-component of the generic initial ideal, with
respect to the reverse lexicographic order, of an ideal generated by
a regular sequence of homogeneous polynomials of degree $d$ is
revlex in a particular, but important, case. Using this property, we
compute the generic initial ideal for several complete intersection
with strong Lefschetz property.

\vspace{5 pt} \noindent \textbf{Keywords:} complete intersection,
generic initial ideal, Lefschetz property.

\vspace{5 pt} \noindent \textbf{2000 Mathematics Subject
Classification:} Primary 13P10, Secondary 13D40,13C40.
\end{abstract}

\begin{center}
\textbf{Introduction.}
\end{center}

Let $K$ be an algebraically closed field of characteristic zero. Let
$S=K[x_{1},\ldots,x_{n}]$ be the polynomial ring in $n$ variables
over $K$. Let $n,d\geq 2$ be two integers. We consider \linebreak
$I=(f_{1},\ldots,f_{n})\subset S$ an ideal generated by a regular
sequence $f_{1},\ldots,f_{n}\in S$ of homogeneous polynomials of
degree $d$. We say that $A=S/I$ is a \emph{$(n,d)$-complete
intersection}. Let $J=Gin(I)$ be the generic initial of $I$, with
respect to the reverse lexicographical (revlex) order (see \cite[\S
15.9]{E}, for details).

We say that a property $(P)$ holds for a generic sequence of
homogeneous polynomials $f_1,f_2,\ldots, f_n \in S$ of given degrees
$d_1,d_2,\ldots,d_n$ if there exists a nonempty open Zariski subset
$U\subset S_{d_1}\times S_{d_2}\times \cdots \times S_{d_n}$ such
that for every $n$-tuple $(f_{1},f_{2},\ldots,f_{n})\in U$ the
property $(P)$ holds. We say that a set of monomials $M\subset S$ is
a {\em revlex set} if, given a monomial $u\in M$, then any other monomial greater than $u$ in revlex order is also in $M$.

For any nonnegative integer $k$, we denote by $J_{k}$ the set of
monomials from $J$ of degree $k$. Conca and Sidman proved that $J_d$
is revlex if $f_{1},\ldots,f_{n}$ is a generic regular sequence,
(see \cite[Theorem 1.2]{CS}). In the first part of this paper, we
prove that $J_{d}$ is a revlex set in another case, namely, when
$f_{i}\in k[x_{i},\ldots,x_{n}]$. It is likely to be true that
$J_{d}$ is revlex for any $(n,d)$-complete intersection, but we do
not have the means to prove this assertion.

We say that a homogeneous polynomial $f$ of degree $s$ is
\emph{semiregular} for $S/I$ if the maps $(S/I)_{t}\stackrel{\cdot
f}{\longrightarrow} (S/I)_{t+s}$ are either injective, either
surjective for all $t\geq 0$. We say that $S/I$ has the \emph{weak
Lefschetz property} (WLP) if there exists a linear form $\ell\in S$,
semiregular on $S/I$. In such case, we say that $\ell$ is a weak
Lefschetz element for $S/I$. A theorem of
Harima-Migliore-Nagel-Watanabe (see \cite{HMNW}) states that $S/I$
has (WLP) in the case $n=3$. We say that $S/I$ has the \emph{strong
Lefschetz property} $(SLP)$ if there exists a linear form $\ell\in
S$ such that $\ell^{b}$ is semiregular on $S/I$ for all integer
$b\geq 1$. In this case, we say that $\ell$ is a strong Lefschetz
element for $S/I$. Harima-Watanabe \cite{HW} and later
Herzog-Popescu \cite{HP}, proved that $S/I$ has (SLP) if $f_{i}\in
k[x_{i},\ldots,x_{n}]$, for all $1\leq i\leq n$.

In the second section of our paper, we compute the generic initial
ideal for some particular cases of $(n,d)$-complete intersections:
$(n=4,d=2)$, $(n=5,d=2)$ and $(n=4,d=3)$. In order to do this, we
suppose in addition that $S/I$ has (SLP). Note that this property
holds for generic complete intersection (see \cite{Par}) and also in
the case when $f_{i}\in k[x_{i},\ldots,x_{n}]$. It was conjectured
that (SLP) holds for any standard complete intersection. A theorem
of Wiebe \cite{W} states that $S/I$ has (WLP) (respectively (SLP))
if and only if $x_n$ is a weak (respectively strong) Lefschetz
element for $S/J$, where $J=Gin(I)$.
As Example $1.9$ show, the hypothesis
$char(K)=0$ and $f_1,\ldots,f_n$ is a regular sequence are
essentials. \vspace{10 pt}

\section{Generic initial ideal for $(n,d)$-complete intersections.}

Let $I=(f_{1},\ldots,f_{n})\subset S = K[x_1,\ldots,x_n]$ be an
ideal generated by a regular sequence $f_{1},\ldots,f_{n}\in S$ of
homogeneous polynomials of degree $d$. Let $J=Gin(I)$ be the generic
initial ideal of $I$, with respect to the revlex order. It is well
known that the Hilbert series of $S/J$ is the same as the Hilbert
series of $S/I$ and moreover, $H(S/J,t) = H(S/I,t) = (1+t + \cdots +
t^{d-1})^{n}$. More precisely, we have:

\begin{prop}
\begin{enumerate}
    \item $H(S/J,k) = \binom{k+n-1}{n-1}$, for $0\leq k\leq d-1$.
    \item $H(S/J,k) = \binom{k+n-1}{n-1} - n\binom{j+n-1}{n-1}$, for $d\leq k\leq \left\lfloor \frac{n(d-1)}{2} \right\rfloor$ and $j=k-d$.
    \item $H(S/J,k) = H(S/J,n(d-1)-k)$, for $k\geq \left\lceil \frac{n(d-1)}{2}\right\rceil$.
\end{enumerate}
\end{prop}

\begin{proof}
Use induction on $n$. Denote $H_{n}(t)=(1+t + \cdots +
t^{d-1})^{n}$. The case $n=1$ is trivial. The induction step follows
from the equality $H_{n}(t) = H_{n-1}(t)(1+t + \cdots + t^{d-1})$.
\end{proof}

\begin{cor}
\begin{enumerate}
  \item $|J_{k}|=0$, for $k\leq d-1$.
  \item $|J_{k}|=n\binom{j+n-1}{n-1}$, for $d\leq k\leq \left\lfloor \frac{n(d-1)}{2} \right\rfloor$ and $j=k-d$.
  \item $|J_{k}|= \binom{\left\lceil \frac{n(d-1)}{2} \right\rceil + j + n - 1}{n-1} -
        \binom{ \left\lfloor\frac{n(d-1)}{2} \right\rfloor -j + n - 1}{n-1} +
        n\binom{\left\lfloor\frac{n(d-1)}{2} \right\rfloor - d - j - n}{n-1}$, for
        $\left\lceil \frac{n(d-1)}{2} \right\rceil \leq k \leq (n-1)(d-1)-1$, where $j=k-\left\lceil \frac{n(d-1)}{2}
        \right\rceil$
  \item $|J_{k}| = \binom{(n-1)d+j}{n-1} - \binom{n-1+d-1-j}{n-1}$, for $(n-1)(d-1)\leq k \leq n(d-1)$, where
        $j=k-(n-1)(d-1)$.
\end{enumerate}
\end{cor}

\begin{proof}
Using $|J_k| = |S_k| - H(S/J,k)$ the proof follows from $1.1$.
\end{proof}

Suppose $f_{i} = \sum_{k=1}^{N}b_{ik}u_{k}$ for $1\leq i \leq n$
where $u_{1},u_{2},\ldots,u_{N}\in S$ are all the monomials of
degree $d$ decreasing ordered in revlex and $N=\binom{d+n-1}{n-1}$.
We denote $u_{k}=x^{\alpha_{k}}$. For example,
$\alpha_1=(d,0,\ldots,0)$, $\alpha_2=(d-1,1,0,\ldots,0)$ etc.

We take a generic transformation of coordinates $x_{i} \mapsto
\sum_{j=1}^{n}c_{ij}x_{j}$ for $i=1,\ldots,n$. Conca and Sidman
proved in \cite{CS} that we may assume that $c_{ij}$ are algebraically
independents over $K$. More precisely, if we consider the field
extension $K\subset L=K(c_{ij}|i,j=\overline{1,n})$ and if we set
\[ F_{i} = f_{i}(\sum_{j=1}^{n}c_{1j}x_{j},\ldots,\sum_{j=1}^{n}c_{nj}x_{j}) \in L[x_{1},\ldots,x_{n}],\;i=1,\ldots,n\]
then $J = Gin(I) = in(F_{1},\ldots,F_{n})\cap S$.

We write $F_{i} = \sum_{j=1}^{n}a_{ij}u_{j} + \cdots $ the monomial
decomposition of $F_{i}$ in $L[x_{1},\ldots,x_{n}]$. With these
notations, we have the following elementary lemma:

\begin{lema}
$J_{d}$ is revlex if and only if the following condition is
fulfilled:
\[ \Delta = \left|
\begin{array}[pos]{ccc}
    a_{11} & \cdots & a_{1n} \\
    \vdots &        & \vdots \\
    a_{n1} & \cdots & a_{nn}
\end{array}
\right|\neq 0.\]
\end{lema}

\begin{proof}
Suppose $\Delta\neq 0$. Since $|J_{d}|=n$, it is enough to show that
$u_{1},\ldots,u_{n}\in J$. Let $A=(a_{ij})_{\tiny
\begin{array}{c} i=1,n\\ j=1,n
\end{array}}$. Since $\Delta=det(A)\neq 0$, $A$ is invertible and we
have
\[A^{-1}
\left( \begin{array}[pos]{c}
    F_1 \\
    \vdots \\
    F_n
\end{array} \right) = \left( \begin{array}[pos]{c}
    H_1 \\
    \vdots \\
    H_n
\end{array} \right),
 \]
 where $H_i = u_i + $ small terms in revlex order. Therefore $LM(H_i)=u_i \in J$, for all $1\leq i\leq n$,
 where $LM(H_i)$ denotes the leading monomial of $H_i$ in the revlex order.

 Conversely, since $u_{1},\ldots,u_{n}\in J_d$, we can find some polynomials $H_i\in L[x_1,\ldots,x_n]$, with
 $LM(H_i)=u_i$, $1\leq i\leq n$, as linear combination of $F_i$'s. 
 If we denote $H_i=\sum_{j=1}^N \widetilde{a}_{ij}u_j$ and $\widetilde{A} = (\widetilde{a}_{ij})_{i,j=1,\ldots,n}$, it
 follows that there exists a map $\psi:L^{n} \rightarrow L^{n}$, given by a matrix $E=(e_{ij})_{i,j=1,\ldots,n}$, such
 that $\widetilde{A} = A\cdot E$. Now, since $det(\widetilde{A})\neq 0$ it follows that $\Delta=det(A)\neq 0$, as
 required.
\end{proof}

\begin{obs}{\em
By the changing of variables $\varphi$ given by $x_{i}\mapsto
\sum_{j=1}^{n}c_{ij}x_{j}$, $x^{\alpha_{k}}$ became
\[ m_{k}:=(\sum_{j=1}^{n}c_{1j}x_{j})^{\alpha_{k1}}\cdots (\sum_{j=1}^{n}c_{nj}x_{j})^{\alpha_{kn}} =
 (\sum_{|t|=\alpha_{k1}}c_{1}^{t}x^{t}) \cdots (\sum_{|t|=\alpha_{kn}}c_{n}^{t}x^{t}) ,\]
where, for any multiindex $t=(t_{1},\ldots,t_{n})$ we denoted
$x^t=x_1^{t_1}\cdots x_n^{t_n}$ and $c_{i}^{t} =
c_{i1}^{t_{1}}\cdots c_{in}^{t_{n}}$. Let $g_{kl}$ be the
coefficient in $c_{ij}$'s of $x^{\alpha_{l}}$ in the monomial
decomposition of $m_{k}$. Using the above writing of $m_k$, we claim
that: \pagebreak
\[(1)\;\; g_{kl} = \tiny \sum_{\begin{array}{c} |t_{1}| = \alpha_{k1},\ldots,|t_{n}| = \alpha_{kn} \\ t_{1}+\cdots+t_{n} = \alpha_{l} \end{array}}
\left[ \binom{\alpha_{k1}}{t_{11}}\cdots \binom{\alpha_{kn}}{t_{n1}}
\right] \left[ \binom{\alpha_{k1}-t_{11}}{t_{12}}\cdots
\binom{\alpha_{kn}-t_{n1}}{t_{n2}}\right] \cdots \] \[ \left[
\binom{\alpha_{k1}-t_{11} - \cdots t_{1n-1}}{t_{1n}}\cdots
\binom{\alpha_{kn}-t_{n1} - \cdots - t_{nn-1}}{t_{nn}} \right] \cdot
c_{1}^{t_{1}}\cdots c_{n}^{t_{n}}. \] \normalsize Indeed, the
monomial $c_{1}^{t_{1}}\cdots c_{n}^{t_{n}}$ appear in the
coefficient of $x^{\alpha_{l}}$ in the expansion of $m_k$ if and
only if $t_1+\cdots+t_n = \alpha_l$ and $|t_{1}| =
\alpha_{k1},\ldots,|t_{n}| = \alpha_{kn}$. Moreover, by Newton
binomial, the coefficient of $x_1^{t_{i1}}\cdots x_n^{t_{in}}$ in
$(\sum_{j=1}^{n}c_{ij}x_{j})^{\alpha_{k1}}$ is
$\binom{\alpha_{k1}}{t_{i1}}
\binom{\alpha_{k1}-t_{i1}}{t_{i2}}\cdots  \binom{\alpha_{k1}-t_{i1}
- \cdots t_{i,n-1}}{t_{in}}c_i^{t_i}$ for any $1\leq i\leq n$, and
thus we proved the claim.

Since $a_{il} = \sum_{k=1}^{N}b_{ik}\cdot g_{kl}$, from the
Cauchy-Binet formula we get:
\[ \Delta =   \sum_{1\leq k_{1} < k_{2} <\ldots < k_{n} \leq N} B_{k_{1},k_{2},\ldots,k_{n}}G_{k_{1},k_{2},\ldots,k_{n}}, \; where \]
\[ B_{k_{1},k_{2},\ldots,k_{n}} =
\left| \begin{array}[pos]{ccc}
    b_{1k_{1}} & \cdots & b_{1k_{n}} \\
    \vdots &  & \vdots \\
    b_{nk_{1}} & \cdots & b_{nk_{n}}
\end{array} \right| \; and \;\;
G_{k_{1},k_{2},\ldots,k_{n}} = \left| \begin{array}[pos]{ccc}
    g_{k_{1}1} & \cdots & g_{k_{n}1} \\
    \vdots &  & \vdots \\
    g_{k_{1}n} & \cdots & g_{k_{n}n}
\end{array} \right|. \]
}\end{obs}

Now, we are able to prove the main result of our paper.

\begin{teor}
If $f_{i}\in K[x_{i},\ldots,x_{n}]$ then $J_{d}$ is revlex. In
particular, if $S/I$ is a monomial complete intersection, then
$J_{d}$ is revlex.
\end{teor}

\begin{proof}
Let $k_{i}=\binom{i+d-1}{d}$, for any $i=1,\ldots,n$. Then
$u_{k_{i}} = x_{i}^{d}$. Recall our notation, $u_k=x^{\alpha_k}$. We
have $b_{11}\neq 0$, otherwise $I = (f_{1},\ldots,f_{n}) \subset
(x_{2},\ldots,x_{n})$ contradicting the fact that $I$ is an Artinian
ideal.  Using a similar argument, we get $b_{ik_{i}}\neq 0$ for all
$1 \leq i \leq n$. Thus, multiplying each $f_{i}$ with
$b_{ik_{i}}^{-1}$, we may assume $b_{ik_{i}}=1$ for all $1 \leq i
\leq n$. In other words, $f_{i}=x_{i}^{d} + f'_{i}$, where $f'_{i}$
contains monomials smaller than $x_{i}^{d}$ in the revlex order.
Also, since $f_{i}\in K[x_{i},\ldots,x_{n}]$ we have $b_{i'k_{i}}=0$
for any $i'>i$. In particular,$B_{k_{1},\ldots,k_{n}}=1$.

In the expansion of the determinant $G_{k_{1},\ldots,k_{n}}$,
appears the term $g_{k_{1}1}\cdot g_{k_{2}2} \cdots g_{k_{n}n} =
r\cdot (c_{11}^{d})(c_{21}^{d-1}c_{22})\cdots (c_{i}^{\alpha_{i}})
\cdots (c_{n}^{\alpha_{n}})$, where $r$ is a nonzero (positive)
integer. Indeed, by $(1)$, we have $g_{11}=c_{11}^d$, $g_{k_{2}2}=d
c_{21}^{d-1}c_{22}$ and, in general, $g_{k_{i}i} = $ some binomial
coefficient $\cdot c_i^{\alpha_i}$. We claim that
$m=(c_{11}^{d})(c_{21}^{d-1}c_{22})\cdots (c_{i}^{\alpha_{i}})
\cdots (c_{n}^{\alpha_{n}})$ doesn't appear again in the expansion
of $\Delta$.

Since $f_{i}\in k[x_{i},\ldots,x_{n}]$, in the monomials in
$(c_{tl})$ of $a_{ij}$ there are no $c_{tl}$'s with $t<i$. Also, all
the monomials of $f'_{i}$ contain variables $x_{t}$ with $t>i$.
Corresponding to them, in $a_{ij}$'s there are $c_{tj}$'s with
$t>i$. Thus in $a_{il}$ the only monomials in $c_{i1},\ldots,c_{in}$
of degree $d$ comes from
$\varphi(x_i^{d})=(\sum_{j=1}^{n}c_{ij}x_j)^d$, the other monomials
being multiples of some $c_{tl}$ with $t>i$. Consequently, in the
expansion of $\Delta$, the monomials of the type
$c_1^{\beta_1}\cdots c_n^{\beta_n}$, where $\beta_1,\ldots,\beta_n$
are multiindices with $|\beta_1|= \cdots = |\beta_n|=d$ comes only
from $\varphi(x_1^d),\ldots,\varphi(x_n^d)$.

On the other hand, for any $1\leq i\leq n$, $c_i^{\alpha_i}$ is unique between the
monomials in $c_{tl}$'s from $\varphi(x_n^d)$, because they are of the type $c_i^{\gamma}$, 
where $\gamma$ is a multiindex with $|\gamma|=d$. 
From these facts, it follows that the
monomial $m$ is unique in the monomial expansion of $\Delta$ and
occurs there with a nonzero coefficient. Thus $\Delta\neq 0$ and by
applying Lemma $1.3$ we complete the proof of the theorem.
\end{proof}

\begin{obs}{\em
In the case $n=2$ and $n=3$, $J_{d}$ is revlex for any
$(n,d)$-complete intersection. Indeed, in the case $n=2$, $J$ itself
is revlex since it is strongly stable. In the case $n=3$, since
$|J_{d}|=3$ and $J$ is strongly stable, it follows that either (a)
$J_{d}=(x_{1}^{d},x_{1}^{d-1}x_{2},x_{1}^{d-2}x_{2}^{2})$, either
(b) $J_{d}=(x_{1}^{d},x_{1}^{d-1}x_{2},x_{1}^{d-1}x_{3})$. But in
the case $(b)$, the map $(S/J)_{d-1} \stackrel{\cdot
x_{3}}{\longrightarrow} (S/J)_{d}$ is not injective, because
$x_{1}^{d-1}\neq 0$ in $(S/J)_{d-1}$ and $x_{1}^{d-1}x_{3}= 0$ in
$(S/J)_{d}$. This is a contradiction with the fact that $x_3$ is a weak
Lefschetz element on $S/J$ and therefore, $J_d$ is
revlex.}\end{obs}

\begin{lema}
(a) $a_{i1}=f_{i}(c_{11},\ldots,c_{n1})$ for all $1\leq i\leq n$.

(b) If $1\leq l\leq n$ is an integer then the sequence
$a_{1l},a_{2l},\ldots,a_{nl}$ is regular as a sequence of
          polynomials in $K[c_{ij}|\;1\leq i,j\leq n]$.
\end{lema}

\begin{proof}
Substituting $x_j=0$ for $j\neq 1$ in $F_i$ we get (a). In order to
prove (b), firstly notice that $a_{11},a_{21},\ldots,a_{n1}$ is a
regular sequence on $K[c_{11},c_{21},\ldots,c_{n1}]$, since $f_1,\ldots,f_n$ is a regular sequence on $K[x_1,\ldots,x_n]$ and $c_{11},c_{21},\ldots,c_{n1}$ are algebraically independent over $K$. 

Let $1\leq l\leq n$ be an integer. We claim that
\[(*)\; \frac{K[c_{ij}|\;1\leq i,j\leq n]}{(a_{1l},\ldots,a_{nl},c_{i1}-c_{ij}\;,\;1\leq i\leq n,\;2\leq j\leq n)} \cong
        \frac{K[c_{11},c_{21},\ldots,c_{n1}]}{(a_{11},a_{21},\ldots,a_{n1})}.\]
Indeed, by $(1)$, if we put $c_{ij}=c_{i1}$ for all $1\leq i\leq n$ and $2\leq j\leq n$ in the expansion of $g_{kl}$
we obtain $r_l \cdot g_{k1}$, where $r_l$ is a strictly positive integer, which depends only on $l$, and therefore, 
$a_{il}$ became $r_l\cdot a_{i1}$. From (*) it follows that $a_{1l},\ldots,a_{nl},c_{i1}-c_{ij}\;for\;1\leq i\leq n,\;2\leq j\leq n$ is a system of parameters for $K[c_{ij}|\;1\leq i,j\leq n]$ and thus $a_{1l},\ldots,a_{nl}$ is a regular sequence on $K[c_{ij}|\;1\leq i,j\leq n]$, so we proved (b).
\end{proof}

As we noticed in Remark $1.6$, for $n=3$, the conclusion of Theorem
$1.5$ holds for any regular sequence $f_1,f_2,f_3$ of homogeneous
polynomials of degree $d$. In the following, we give another proof
of this, without using the fact that $S/(f_1,f_2,f_3)$ has the
(WLP), i.e. $x_3$ is a weak Lefschetz element for $S/J$. Also, we
get the same conclusion for the case $n=4$ and $d=2$. However, this
approach do not works in the general case.

\begin{prop}
(a) If $f_1,f_2,f_3\in K[x_1,x_2,x_3]$ is a regular sequence of
homogeneous polynomials of degree $d\geq 2$, $I=(f_1,f_2,f_3)$
    and $J=Gin(I)$, the generic initial ideal of $I$, with respect to the reverse lexicographical order, then
    $J_d$ is a revlex set.

(b) If $f_1,f_2,f_3,f_4\in K[x_1,x_2,x_3,x_4]$ is a regular sequence
of homogeneous polynomials of degree $2$, $I=(f_1,f_2,f_3,f_4)$
    and $J=Gin(I)$, the generic initial ideal of $I$, with respect to the reverse lexicographical order, then
    $J_2$ is a revlex set.
\end{prop}

\begin{proof}
(a) Let $A=(a_{ij})_{i,j=\overline{1,3}}$. Since $Gin(f_1,f_2)$ is
strongly stable, it follows by Lemma $1.3$ that
    $\Delta_{3}=\left| \begin{array}[pos]{cc} a_{11} & a_{12}\\ a_{21} & a_{22} \end{array} \right| \neq 0$.
    Analogously, $\Delta_{2}=\left| \begin{array}[pos]{cc} a_{11} & a_{12}\\ a_{31} & a_{32} \end{array} \right| \neq 0$ and $\Delta_{1}=\left| \begin{array}[pos]{cc} a_{21} & a_{22}\\ a_{31} & a_{32} \end{array} \right| \neq 0$.
    We have $\Delta = a_{13}\Delta_{1} - a_{23}\Delta_{2} + a_{33}\Delta_{3}$.
    Suppose $\Delta = 0$. It follows $a_{13}\Delta_{1} = a_{23}\Delta_{2} - a_{33}\Delta_{3}$
    and therefore, since $a_{13},a_{23},a_{33}$ is a regular sequence in $K[c_{ij}| i,j=\overline{1,3}]$,
    we get $\Delta_{1} \in (a_{23},a_{33})$. The first three monomials of degree $d$ in revlex order are $x_1^d$,
    $x_1^{d-1}x_2$ and $x_1^{d-2}x_2^2$. It follows that the degree of
    $a_{i1}$, $a_{i2}$ and $a_{i3}$ in $c_{21},c_{22},c_{23}$ is $0$, $1$, respectively $2$,
    for any $1\leq i\leq 3$. Therefore, the degree of $\Delta_{1}$ in
    the variables $c_{21},c_{22},c_{23}$ is $1$, but the degree of $a_{23}$ and $a_{33}$ in $c_{21},c_{22},c_{23}$ is
    $2$, which is impossible, since $\Delta_{1} \in (a_{23},a_{33})$.

(b) Let $A=(a_{ij})_{i,j=\overline{1,4}}$. Since any three
polynomials from $f_1,f_2,f_3,f_4$ form a regular sequence, it
follows from (a) that any $3\times 3$
    minor of the matrix $\widetilde{A}=(a_{ij})_{\tiny \begin{array}{c} i=\overline{1,4} \\ j=\overline{1,3}
    \end{array} \normalsize}$ is
    nonzero. Let $\Delta_i$ be the minor obtained from $\widetilde{A}$ by erasing the $i$-row. Suppose $\Delta=0$. It
    follows that $a_{14}\Delta_1 = a_{24}\Delta_2 - a_{34}\Delta_3 + a_{44}\Delta_4$ and therefore, since
    $a_{14},a_{24},a_{34},a_{44}$ is a regular sequence in $K[c_{ij}| i,j=\overline{1,4}]$, we get
    $\Delta_1 \in (a_{24},a_{34},a_{44})$. Since the first $4$ monomials in
    revlex are $x_1^2,x_1x_2,x_2^2,x_1x_3$, we get a contradiction from the fact that
    the degree of $\Delta_1$ in the variables $c_{31},c_{32},c_{33},c_{34}$ is zero, but the degree of
    $a_{24},a_{34},a_{44}$ in $c_{31},c_{32},c_{33},c_{34}$ is $1$.
\end{proof}


\begin{obs}{\em
The hypothesis that $K$ is a field with $char(K)=0$ is essential.
Indeed, suppose $char(K)=p$ and $I = (x_1^{p},x_2^{p})\subset
K[x_1,x_2]$. Then, simply using the definition of the generic
initial ideal, we get $Gin(I)=I$ and, obviously,
$I_{p}=\{x_1^{p},x_2^{p}\}$ is not revlex.

Also, the hypothesis that $f_1,\ldots,f_n$ is a regular sequence of
homogeneous polynomials is essential. Let $I=(f_1,f_2,f_3)\subset
K[x_1,x_2,x_3]$, where $f_1=x_1^2$, $f_2=x_1x_2$ and $f_3=x_1x_3$.
In order to compute the generic initial ideal of $I$ we can take a
generic transformation of coordinates with an upper triangular
matrix, i.e. $x_1 \mapsto x_1,\; x_2 \mapsto x_2+c_{12}x_1,\;
x_3\mapsto x_3 + c_{23}x_2 + c_{13}x_1$, where $c_{ij}\in K$ for all
$i,j$ (see \cite[\S 15.9]{E}). We get
\[ F_{1}(x_1,x_2,x_3):=f_{1}(x_1,x_2+c_{12}x_1,x_3 + c_{23}x_2 + c_{13}x_1)= x_1^{2}, \]
\[ F_{2}(x_1,x_2,x_3):=f_{2}(x_1,x_2+c_{12}x_1,x_3 + c_{23}x_2 + c_{13}x_1)= c_{12}x_1^{2}+x_1x_2, \]
\[ F_{3}(x_1,x_2,x_3):=f_{3}(x_1,x_2+c_{12}x_1,x_3 + c_{23}x_2 + c_{13}x_1)= c_{13}x_1^{2}+c_{23}x_1x_2 + x_1x_3. \]
The generic initial ideal of $I$, $J=in(F_1,F_2,F_3)$ satisfies $J_2=I_{2}$, but $I_2$ is not revlex.}
\end{obs}

\section{Several examples of computation of the Gin.}

Let $I=(f_{1},\ldots,f_{n})\subset S = K[x_1,\ldots,x_n]$ be an
ideal generated by a regular sequence $f_{1},\ldots,f_{n}\in S$ of
homogeneous polynomials of degree $d$. Let $J=Gin(I)$ be the generic
initial ideal of $I$, with respect to the revlex order.

In \cite{C}, the case $n=3$ and $d\geq 2$ is treated completely,
when $S/(f_1,f_2,f_3)$ has (SLP). More precisely, if $d$ is odd,
then
\[ J = ( x_{1}^{d-2}\{x_{1},x_{2}\}^{2}, x_{1}^{d-2j-1}x_{2}^{3j+1}, x_{1}^{d-2j-2}x_{2}^{3j+2}\; for\; 1\leq j \leq  \frac{d-3}{2} , \; x_{2}^{\frac{3d-1}{2}}, x_{3}x_{2}^{\frac{3d-3}{3}}, \]\[
x_{3}^{2j+1}x_{1}^{2j}x_{2}^{\frac{3d-3}{2}-3j}, \ldots,
x_{3}^{2j+1}x_{2}^{\frac{3d-3}{2}-j} , 1\leq j \leq \frac{d-3}{2} \;
 ,x_{3}^{d-2+2j}\{x_{1},x_{2}\}^{d-j}, 1\leq j\leq d) \]
\[ or\;\; J = ( x_{1}^{d-2}\{x_{1},x_{2}\}^{2}, x_{1}^{d-2j-1}x_{2}^{3j+1}, x_{1}^{d-2j-2}x_{2}^{3j+2}\; for\; 1\leq j \leq  \frac{d-4}{2}, \; x_{1}x_{2}^{\frac{3d-4}{2}}, x_{2}^{\frac{3d-2}{2}},\; \] \[
x_{3}^{2j}x_{1}^{2j-1}x_{2}^{\frac{3d}{2}-3j}, \ldots,
x_{3}^{2j}x_{2}^{\frac{3d-2}{2}-j},  1\leq j\leq\frac{d-2}{2}
 ,x_{3}^{d-2+2j}\{x_{1},x_{2}\}^{d-j}, 1\leq j\leq d) \]
if $d$ is even (see \cite[Proposition 3.3]{C}).

In the following, we discuss some particular cases with $n\geq 4$.

\paragraph{The case $n=4$, $d=2$.}

We assume that $S/I$ has (SLP). From Wiebe's Theorem, it follows
that $x_{4}$ is a strong Lefschetz element for $S/J$. For a positive
integer $k$, we denote $Shad(J_{k}) = \{x_1,\ldots,x_n\}J_{k}$. We
have $H(S/J,t) = (1+t)^{4} = 1+4t+6t^{2}+4t^{3}+t^{4}$.

We have $|J_{2}|=4$. From Proposition $1.8$, $J_2$ is revlex,
therefore
\[ J_{2}=\{x_{1}^{2},x_{1}x_{2},x_{2}^{2},x_{1}x_{3}\} = \{\{x_{1},x_{2}\}^{2},x_{1}x_{3}\} .\]
We have $|Shad(J_{2})| = 12$. On the other hand, $|J_{3}|=16$, so we
need to add $4$ new generators at $Shad(J_{2})$ to get $J_{3}$. If
we add a new monomial which is divisible by $x_{4}^{2}$, then the
map $(S/J)_{1} \stackrel{\cdot x_{4}^{2}}{\longrightarrow}
(S/J)_{3}$, will be no longer injective. Since $|(S/J)_{1}| =
|(S/J)_{3}|$, we get a contradiction with the fact that $x_{4}$ is a
strong Lefschetz element for $S/J$. But there exists only $16$
monomials in $S$ which are not multiple of $x_4^2$. Thus
\[J_{3} = \{\{x_{1},x_{2},x_{3}\}^{3}, x_{4}\{x_{1},x_{2},x_{3}\}^{2}\},\;and\; therefore\]
\[ Shad(J_{3})=\{\{x_{1},x_{2},x_{3}\}^{4}, x_{4}\{x_{1},x_{2},x_{3}\}^{3}, x_{4}^{2}\{x_{1},x_{2},x_{3}\}^{2}\}. \]
Since $|Shad(J_{3})| = 31$ and $|J_{4}| = |S_{4}|-|(S/J)_{4}| =35 -
1 = 34$ we have to add $3$ new generators at $Shad(J_{3})$ in order
to get $J_{4}$. Since $J$ is strongly stable, these new generators
are $x_{4}^{3}x_{1}$, $x_{4}^{3}x_{2}$ and $x_{4}^{3}x_{3}$. So
\[ J_{4}=\{x_{1},x_{2},x_{3},x_{4}\}^{4}\setminus \{x_{4}^{4}\}.\;We\;get\;\;
   Shad(J_{4})= \{x_{1},x_{2},x_{3},x_{4}\}^{5}\setminus \{x_{4}^{5}\}\]
   and since $J_{5}=S_{5}$ it follows that we must add $x_{4}^{5}$ at $Shad(J_{4})$ to obtain $J_{5}$. From now one,
   we cannot add any new monomial. $J$ is
   the ideal generated by all monomials added at some step $k$ to $Shad(J_{k})$, thus
   we proved the following proposition:

\begin{prop}
If $I=(f_{1},f_{2},f_{3},f_{4})$ is an ideal generated by a regular
sequence of homogeneous polynomials $f_{1},f_{2},f_{3},f_{4} \in
S=k[x_{1},x_{2},x_{3},x_{4}]$ of degree $2$ such that the algebra
$S/I$ has (SLP) then the generic initial ideal of $I$ with respect
to the revlex order is
\[J = (x_{1}^{2},\; x_{1}x_{2},\; x_{2}^{2},\; x_1x_3,\;
       x_{2}x_{3}^{3},\; x_{3}^{3}, \; x_{3}^{2}x_{4}, \; x_{3}^{2}x_{4},\;
       x_{4}^{3}x_{1},\; x_{4}^{3}x_{2},\; x_{4}^{3}x_{3},\; x_{4}^{5}). \]
In particular, this assertion holds for a generic sequence of
homogeneous polynomials $f_{1},f_{2},f_{3},f_{4}\in S$ or if
$f_{i}\in k[x_{i},\ldots,x_{4}]$, $1\leq i\leq 4$.
\end{prop}

\paragraph{The case $n=5$, $d=2$.}

In the following, we suppose that $S/I$ has (SLP), so $x_{5}$ is a
strong Lefschetz element for $S/J$. Also, we suppose that $J_{2}$ is
revlex. We have $H(S/J,t) = (1+t)^{5} = 1+5t+10t^{2}+10t^{3}+5t^{4}
+ t^{5}$. We have $|J_{2}|=5$. Since $J_{2}$ is revlex from the
assumption, we have $J_{2} =
\{\{x_1,x_2\}^{2},x_{3}\{x_{1},x_{2}\}\}$. So \[ Shad(J_{2}) = \{
\{x_{1},x_{2}\}^{3}, \{x_{1},x_{2}\}^{2}\{x_3,x_4,x_5\},
x_3\{x_1,x_2\}\{x_3,x_4,x_5\}\}.\] We have $|Shad(J_{2})|=19$. On
the other hand $|J_{3}| = |S_{3}| - |(S/J)_{3}| = 35 - 10 = 25$, so
we must add $6$ new generators, from a list of $16$ monomials, at
$Shad(J_{2})$ to get $J_{3}$.

Since $x_{5}$ is a strong Lefschetz element for $S/J$ it follows
that we cannot add any monomial of the form $x_{5}\cdot m$, where
$m$ is nonzero in $(S/J)_{2}$ because, in that case, the map
$(S/J)_{2} \stackrel{\cdot x_{5}}{\rightarrow} (S/J)_{3}$ will be no
longer injective. But there are $|(S/J)_{2}| = 10$ such monomials
$m$. Therefore, we must add the remaining $6$ monomials,
$x_{3}^{3},x_{3}^{2}x_{4}, x_{1}x_{4}^{2},x_{2}x_{4}^{2},
x_{3}x_{4}^{2}, x_{4}^{3}$. Thus
\[ J_{3} = \{\{x_{1},x_{2},x_{3},x_{4}\}^{3} , x_{5}(\{x_{1},x_{2},x_{3}\}^{2}\setminus\{x_{3}^{2}\})\}.\;Therefore:\]
\[ Shad(J_{3}) =\{ \{x_{1},x_{2},x_{3},x_{4}\}^{4} , x_{5}\{x_{1},x_{2},x_{3},x_{4}\}^{3} ,
x_{5}^{2}( \{x_{1},x_{2},x_{3}\}^{2}\setminus\{x_{3}^{2}\})\}. \] We
have $|Shad(J_{3})|=60$ and $|J_{4}| = |S_{4}| - |(S/J)_{4}| = 70 -
5 = 65$. So we need to add $5$ new generators at $Shad(J_{3})$ to
get $J_{4}$. If we add a monomial which is divisible by $x_{5}^{3}$
we obtain a contradiction from the fact that the map $(S/J)_{1}
\stackrel{\cdot x_{5}^{3}}{\rightarrow} (S/J)_{4}$ is no longer
injective. Therefore, we must add: $x_{3}^{2}x_{5}^{2},
x_{1}x_{4}x_{5}^{2}, x_{2}x_{4}x_{5}^{2}, x_{3}x_{4}x_{5}^{2},
x_{4}^{2}x_{5}^{2}$, and so
\[ J_{4} = \{ \{x_{1},x_{2},x_{3},x_{4}\}^{4} , x_{5}\{x_{1},x_{2},x_{3},x_{4}\}^{3} ,
x_{5}^{2} \{x_{1},x_{2},x_{3},x_{4}\}^{2} \}. \]\[ So
\;\;Shad(J_{4}) = \{ \{x_{1},x_{2},x_{3},x_{4}\}^{5} , \cdots ,
x_{5}^{3} \{x_{1},x_{2},x_{3},x_{4}\}^{2} \}.\] 

We have $|J_{5}| - |Shad(J_{4})| = 4$, so we must add $4$ new generators at $Shad(J_{4})$ to
get $J_{5}$. Since $J$ is strongly stable, these new generators are:
$x_{5}^{4}x_{1},x_{5}^{4}x_{2},x_{5}^{4}x_{3},x_{5}^{4}x_{4}$. Therefore $J_{5} = \{ \{x_{1},\ldots,x_{5}\}^{5}\setminus\{x_{5}^{5}\}\}$. Finally, we must add $x_5^6$ to $Shad(J_5)$ in
order to obtain $J_6$. We proved the following proposition, with the help of 
\cite[Theorem 1.2]{CS} and Theorem $1.5$.

\begin{prop}
If $I=(f_{1},f_{2},\ldots,f_{5})\subset K[x_1,\ldots,x_5]$ is an
ideal generated by a generic (regular) sequence of homogeneous
polynomials of degree $2$ or if $f_{1},f_{2},\ldots,f_{5}$ is a
regular sequence of homogeneous polynomials of degree $2$ with
$f_{i}\in K[x_{i},\ldots,x_{5}]$ for $i=1,\ldots,5$ then $J=Gin(I)$
the generic initial ideal of $I$ with respect to the revlex order
is:
\[ J = (x_{1}^{2},\; x_{1}x_{2},\; x_{2}^{2},\;x_{1}x_{3},\; x_{2}x_{3}, \;
       x_{3}^{3},\; x_{3}^{2}x_{4},\;  x_{1}x_{4}^{2},\; x_{2}x_{4}^{2},\;  x_{3}x_{4}^{2},\;  x_{4}^{3},\;\]\[
       x_{3}^{2}x_{5}^{2},\;x_{1}x_{4}x_{5}^{2},\; x_{2}x_{4}x_{5}^{2},\;x_{3}x_{4}x_{5}^{2},\;  x_{4}^{2}x_{5}^{2},\;
       x_{5}^{4}x_{1},\; x_{5}^{4}x_{2},\; x_{5}^{4}x_{3},\; x_{5}^{4}x_{4},\;  x_{5}^{6} ) \]
\end{prop}

\paragraph{The case $n=4$, $d=3$.}

We suppose that $S/I$ has (SLP), so $x_{4}$ is a strong Lefschetz
element for $S/J$. Also, we suppose that $J_{3}$ is revlex.

We have $H(S/J,t) = (1+t+t^{2})^{4} = (1+2t+3t^{2}+2t{3}+t^{4})^{2}
= $
\[ = 1 +4t + 10t^{2} +16t^{3} + 19t^{4} + 16t^{5} + 10t^{6} + 4t^{7} +t^{8}.\]
Since $|J_{3}|=4$ and $J_{3}$ is revlex, it follows that $J_{3} =  \{x_{1},x_{2}\}^{3}$.
Therefore, we have \linebreak $Shad(J_{3}) = \{ \{x_{1},x_{2}\}^{4}, \{x_{1},x_{2}\}^{3}\{x_{3},x_{4}\}\}$. 
Since $|J_{4}| - |Shad(J_{3})|= 4$, we must add $4$ new generators to $Shad(J_{3})$ to
obtain $J_{4}$. Since $x_{4}$ is a strong Lefschetz element for
$S/J$ we cannot add any monomial of the form $x_{4}\cdot m$, where
$m\neq 0$ in $J_{3}$. Therefore, since $J$ is strongly stable, we
have to choose $3$ monomials from the list
$x_{3}^{2}\{x_{1},x_{2}\}^{2},x_{3}^{3}\{x_{1},x_{2}\},x_{3}^{4}$. There
are two different chooses: either we add (I) $x_{3}^{2}\{x_{1},x_{2}\}^{2}$, either (II) $x_{3}^{2}x_{1}\{x_{1},x_{2},x_{3}\}$.

In the case (I), we get $J_{4}=\{ \{x_{1},x_{2}\}^{4},
\{x_{1},x_{2}\}^{3}\{x_{3},x_{4}\}, x_{3}^{2}\{x_{1},x_{2}\}^{2}\}$, so
\[Shad(J_{4}) = \{ \{x_{1},x_{2}\}^{5}, \{x_{1},x_{2}\}^{4}\{x_{3},x_{4}\}, \{x_{1},x_{2}\}^{3}\{x_{3},x_{4}\}^{2}, x_{3}^{2}\{x_{3},x_{4}\}\{x_{1},x_{2}\}^{2}\}. \]
Since $|J_{5}|-|Shad(J_{4})|=40-34=6$, we need to add $6$ new generators
at $Shad(J_{4})$ to get $J_{5}$. Since $x_{4}$ is a strong Lefschetz
element for $S/J$ we cannot add any monomial of the form
$x_{4}^{2}m$, where $m$ is a nonzero monomial in $J_{3}$. So, we
must add: $x_{3}^{4}\{x_{1},x_{2},x_{3}\},
x_{4}x_{3}^{3}\{x_{1},x_{2},x_{3}\}$. Thus
$J_{5} = \{ \{ x_{1},x_{2},x_{3}\}^{5}, x_{4}\{ x_{1},x_{2},x_{3}\}^{4}, x_{4}^{2}\{x_{1},x_{2}\}^{3}\}$. 

In the case (II), we have $J_{4}=\{ \{x_{1},x_{2}\}^{4},
\{x_{1},x_{2}\}^{3}\{x_{3},x_{4}\},
x_{1}x_{3}^{2}\{x_{1},x_{2},x_{3}\} \}$, so $Shad(J_{4})$ is the set
$\{ \{x_{1},x_{2}\}^{5}, \{x_{1},x_{2}\}^{4}\{x_{3},x_{4}\}, \{x_{1},x_{2}\}^{3}\{x_{3},x_{4}\}^{2}, x_{3}^{2}x_{1}\{x_{3},x_{4}\}\{x_{1},x_{2}\}, x_{3}^{3}x_{1}\{x_{3},x_{4}\}\}$.
Since $|J_{5}|-|Shad(J_{4})|=40-34=6$, we must add $6$ new generators at $Shad(J_{4})$ to get $J_{5}$. 
Since $x_{4}$ is a strong-Lefschetz element for $S/J$, we cannot add any monomial of
the form $x_{4}^{2}m$, where $m\neq 0$ in $J_{3}$. So, we must add:
$x_{3}^{3}x_{2}^{2}, x_{3}^{4}x_{2}, x_{3}^{5},
x_{4}x_{3}^{2}x_{2}^{2}, x_{4}x_{3}^{3}x_{2}, x_{4}x_{3}^{4}$. Thus
\[ J_{5} = \{ \{ x_{1},x_{2},x_{3}\}^{5}, x_{4}\{ x_{1},x_{2},x_{3}\}^{4}, x_{4}^{2}\{x_{1},x_{2}\}^{3}\}, \]
the same as in the case (I). Thus, in both cases (I) and (II), we get:
\[Shad(J_{5}) = \{ \{ x_{1},x_{2},x_{3}\}^{6}, x_{4}\{ x_{1},x_{2},x_{3}\}^{5}, x_{4}^{2}\{x_{1},x_{2},x_{3}\}^{4} , x_{4}^{3}\{x_{1},x_{2}\}^{3}\}. \]
Since $|Shad(J_{5})| = |S_{6}|-16$ and $|J_{6}| = |S_{6}| - 10$, we must add $6$ new generators to
$Shad(J_{5})$ in order to obtain $J_6$. Since $x_{4}$ is a strong-Lefschetz element for $S/J$, these new
generators are not divisible by $x_{4}^{4}$. So, we add
$x_{4}^{3}x_{3}\{x_{1},x_{2}\}^{2},x_{4}^{3}x_{3}^{2}\{x_{1},x_{2}\},
x_{4}^{3}x_{3}^{3}$ and thus,
\[ J_{6} = \{ \{x_{1},x_{2},x_{3}\}^{6}, x_{4}\{ x_{1},x_{2},x_{3}\}^{5}, x_{4}^{2}\{x_{1},x_{2},x_{3}\}^{4},x_{4}^{3}\{x_{1},x_{2},x_{3}\}^{3}\}. \;So\]
\[ Shad(J_{6}) = \{ \{x_{1},x_{2},x_{3}\}^{7}, x_{4}\{ x_{1},x_{2},x_{3}\}^{6}, \ldots, x_{4}^{4}\{x_{1},x_{2},x_{3}\}^{3}\}. \]
$|S_{7}|-|Shad(J_{6})| = 6+4 = 10$ and $|S_{7}| - |J_{7}|=4$, so we
must add $6$ new generators at $Shad(J_{6})$ to get $J_{7}$. Using
the same argument, these new generators must be
$x_{4}^{5}\{x_{1},x_{2},x_{3}\}^{2}$ and therefore
$J_{7} = \{ \{x_{1},x_{2},x_{3}\}^{7}, x_{4}\{ x_{1},x_{2},x_{3}\}^{6}, \ldots, x_{4}^{5}\{x_{1},x_{2},x_{3}\}^{2} \}$. We get
\[ Shad(J_{7}) = \{ \{x_{1},x_{2},x_{3}\}^{8}, x_{4}\{ x_{1},x_{2},x_{3}\}^{7}, \ldots, x_{4}^{6}\{x_{1},x_{2},x_{3}\}^{2}\}. \]
Since $|S_{8}| - |Shad(J_{7})| = 4$ and $|S_{8}|-|J_{8}| = 1$, we
must add $3$ new generators at $Shad(J_{7})$ in order to get $J_8$.
Since $x_{4}$ is strong-Lefschetz, these new generators are $x_{4}^{7}\{x_{1},x_{2},x_{3}\}$, so
$J_8 = \{x_{1},x_{2},x_{3},x_{4}\}^{8}\setminus\{x_{4}^{8}\}$. Finally, we must add $x_{4}^{9}$ to
$Shad(J_8)$ in order to obtain $J_9$. We proved the following proposition, with the help of 
\cite[Theorem 1.2]{CS} and Theorem $1.5$.

\begin{prop}
If $I=(f_{1},f_{2},f_{3},f_{4})\subset K[x_1,x_2,x_3,x_4]$ is an
ideal generated by a generic (regular) sequence of homogeneous
polynomials of degree $3$ or if $f_{1},f_{2},f_{3},f_{4}$ is a
regular sequence of homogeneous polynomials of degree $3$ with
$f_{i}\in k[x_{i},\ldots,x_{4}]$, for $i=1,\ldots,4$, then
$J=Gin(I)$ the generic initial ideal of $I$ with respect to the
revlex order has one of the following forms:
\[(I)\;\;\;\;\;\;\;\;\;\;\;\;\;\;\;\;\;\;\;\; J = (\{x_{1},x_{2}\}^{3},\;x_{3}^{2}\{x_{1},x_{2}\}^{2},\; x_{3}^{4}\{x_{1},x_{2},x_{3}\},\;
                  x_{4}x_{3}^{3}\{x_{1},x_{2},x_{3}\}, \]
          \[x_{4}^{3}x_{3}\{x_{1},x_{2}\}^{2},x_{4}^{3}x_{3}^{2}\{x_{1},x_{2}\}, x_{4}^{3}x_{3}^{3},\;
            x_{4}^{5}\{x_{1},x_{2},x_{3}\}^{2},\; x_{4}^{7}\{x_{1},x_{2},x_{3}\},\; x_{4}^{9} ) \]
\[(II)\; J = (\{x_{1},x_{2}\}^{3},\;x_{3}^{2}x_{1}\{x_{1},x_{2},x_{3}\},\;x_{3}^{3}x_{2}^{2},\; x_{3}^{4}x_{2},\;
                   x_{3}^{5},\; x_{4}x_{3}^{2}x_{2}^{2},\; x_{4}x_{3}^{3}x_{2},\; x_{4}x_{3}^{4},  \]
          \[x_{4}^{3}x_{3}\{x_{1},x_{2}\}^{2},x_{4}^{3}x_{3}^{2}\{x_{1},x_{2}\}, x_{4}^{3}x_{3}^{3},\;
            x_{4}^{5}\{x_{1},x_{2},x_{3}\}^{2},\;x_{4}^{7}\{x_{1},x_{2},x_{3}\}, \; x_{4}^{9}) \]
\end{prop}

\begin{obs}
It seems Conca-Herzog-Hibi noticed in \cite{CHH}, page $838$, that, if $f_{1},f_{2},f_{3},f_{4}$
is a generic sequence of homogeneous polynomials of degree $3$ then the generic initial ideal $J$
has the form (I), and $J=Gin(x_1^3,x_2^3,x_3^3,x_4^3)$ has the form (II).
\end{obs}

\vspace{2mm} \noindent {\footnotesize
\begin{minipage}[b]{10cm}
 Mircea Cimpoeas, Junior Researcher\\
 Institute of Mathematics of the Romanian Academy\\
 Bucharest, Romania\\
 E-mail: mircea.cimpoeas@imar.ro
\end{minipage}}
\end{document}